\documentclass[11pt]{articlefederico}
\usepackage{amsthm}
\usepackage{amsmath}
\usepackage{amssymb}
\usepackage{graphicx}
\usepackage{graphics}

\newtheorem{theorem}{Theorem}[section]
\newtheorem{lemma}[theorem]{Lemma}

\newtheorem{definition}[theorem]{Definition}
\newtheorem{proposition}[theorem]{Proposition}

\newtheorem{example}[theorem]{Example}

\newcommand{\A}{{\cal A}}
\newcommand{\B}{{\cal B}}

\newcommand{\swap}{\bf{swap}}

\begin{document}
\title{\textsf{When do two planted graphs have\\ the same cotransversal matroid?}}
\author{\textsf{Federico Ardila\footnote{Department of Mathematics, San Francisco State University, \textsf{federico@math.sfsu.edu}. \newline Supported by National Science Foundation grant DMS-0801075.} \qquad Amanda Ruiz\footnote{Department of Mathematics, San Francisco State University, \textsf{alruiz@sfsu.edu}. \newline Supported by an Achievement Rewards for College Scientists (ARCS) Scholarship.}}}
\date{}
\maketitle

\begin{abstract} 
Cotransversal matroids are a family of matroids that arise from planted graphs. 
We prove that two planted graphs give the same cotransversal matroid if and only if they can be obtained from each other by a series of local moves.
  \end{abstract}


\section{\textsf{Introduction}}
\label{intro}

Cotransversal matroids are a family of matroids that arise from planted graphs. The goal of this short note is to describe when two planted graphs give rise to the same cotransversal matroid. 

The paper is organized as follows. In Section \ref{sec:prelim} we recall some basic definitions and facts in matroid theory, including the notions of cotransversal and transversal matroids. In Sections \ref{swap} and \ref{sec:max} we introduce the operations of \emph{swapping} and \emph{saturating} on a planted graph, and prove that they preserve the cotransversal matroid. (Theorems \ref{preserve} and \ref{unique}) In Section \ref{sec:dual} we prove a crucial lemma on transversal matroids. 
Finally in Section \ref{proof} we prove our main result: two planted graphs give rise to the same cotransversal matroid if and only if their saturations can be obtained from each other by a series of swaps. (Theorem \ref{iso}) 

This paper is inspired by and analogous to Whitney's work on presentations of graphical matroids. He showed \cite{whitney} that two graphs give rise to the same graphical matroid if and only if they can be obtained from each other by repeatedly applying three operations. Our main theorem is also analogous to Bondy \cite{bondy} and Mason's \cite{mason} elegant theorem that a transversal matroid has a unique maximal presentation. In Sections \ref{sec:max} and \ref{sec:dual} we will explain how our theorem and theirs are connected by matroid duality, and we will see the need to resolve several subtleties that do not arise in that dual setting.




\section{\textsf{Preliminaries}}\label{sec:prelim}

Matroids can be thought of as a notion of independence, which generalizes various notions of independence occuring in linear algebra, field theory, graph theory, matching theory, among others. We begin by recalling some basic notions of the theory of matroids. For a more thorough introduction, we refer the reader to \cite{ardilavideos, oxley, white}.

\begin{definition}
A \emph{matroid} $(E, \B)$ consists of a finite set $E$ and a nonempty family $\B$ of subsets of $E$, called \emph{bases}, with the following property:
If $B_a,B_b\in\mathcal{B}$ and $x\in B_a-B_b$, then there exists $y\in B_b-B_a$ such that $(B_a-x)\cup y\in \mathcal{B}.$
\end{definition}

A prototypical example of a matroid consists of a finite collection of vectors $E$ spanning a vector space $V$, and the collection $\B$ of subsets of $E$ which are bases of $V$.

Matroids have a useful notion of duality, as follows.

\begin{definition}
If $M=(E,\B)$ is a matroid then $\B^\ast=\{E-B\mid B\in \B\}$ is also the collection of bases for a matroid $M^\ast=(E,\B^\ast)$, called the \emph{dual} of M. 
\end{definition}
 
Notice that $(M^*)^*=M$. This allows us to talk about \emph{pairs of dual matroids}. 

Duality behaves beautifully with respect to many of the natural concepts on matroids. In particular, the general theory makes it straighforward to translate many notions and results (\emph{e.g.} definitions, constructions, and theorems) about $M$ into ``dual" notions and results about $M^*$.

\subsection{\textsf{Cotransversal and transversal matroids}}

We are particularly interested in two families of matroids arising in graph theory and matching theory. First we define \emph{cotransversal matroids}, which are the main object of study of this paper.
A vertex of a directed graph $G$ is called a \emph{sink} if it has no outgoing edges. A \emph{routing} is a set of vertex-disjoint directed paths in $G$.

\begin{definition}
A \emph{planted graph} $(G,B)$ is a directed graph $G$ with vertex set $V$ having no loops or parallel edges, together with a specified set of sinks $B\subseteq V$. 
\end{definition}

\begin{theorem} \cite{mason2, oxley}
Given a planted graph $(G,B)$ on $V$, there is a matroid $L(G,B)$ on $V$ whose bases are 
the sets of $|B|$ vertices that can be routed to $B$ through vertex-disjoint directed paths.   
\end{theorem}

Any matroid $M$ that arises in this way is called \emph{cotransversal}, and a planted graph giving rise to it is called a \emph{presentation} of $M$.

\begin{figure}[h]
\centering
\includegraphics[scale=.65]{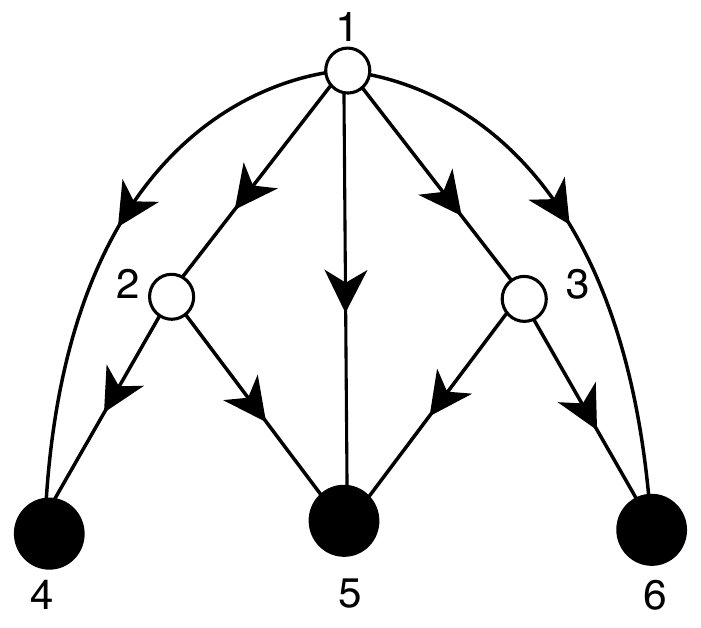}
\caption{A planted graph $(G,B)$ with  $B=\{4,5,6\}$.
\label{move5}}
\end{figure}

\begin{example} \label{ex1} Figure \ref{move5} shows a planted graph $G$ with a specified set of sinks, $B=\{4,5,6\}$. The bases of the cotransversal matroid $M=L(G,B)$ are all $3$-subsets of $\{1,2,3,4,5,6\}$ except $245$ and $356$.
\end{example}

Now we define \emph{transversal matroids}, another important family.

\begin{definition}\label{trans} Let  $S$ be a finite set. Let $\A=\{A_1,\ldots, A_r\}$ be a family
of subsets of $S$.    A \emph{system of distinct representatives} \emph{(SDR)} of $\A$ is a choice of an element $a_i \in A_i$ for each $i$ such that $a_i \neq a_j$ for $i \neq j$. 
A \emph{transversal} is a set which can be ordered to obtain an SDR. 
\end{definition}

\begin{theorem} \cite{oxley}
Given a family  $\A=\{A_1,\ldots, A_r\}$ of subsets of $S$, there is a matroid on $S$ whose bases are the transversals of $\A$. 
\end{theorem}

A matroid that arises in this way is called a \emph{transversal matroid}, and $\A$ is called a \emph{presentation} of it. We can also view $\A = \{A_1, \ldots, A_r\}$ as a bipartite graph between the ``top" vertex set $[r]=\{1,\ldots,r\}$ and the ``bottom" vertex set $S$, where top vertex $i$ is connected to the elements of $A_i$ for $1 \leq i \leq r$. The SDRs of $\A$ become maximal
matchings of $[r]$ into $S$ in this bipartite graph. We will use these two points of view interchangeably.

\begin{example} \label{ex2}
Let $S =\{1, \ldots, 6\}$ and $\A=\{\{1, 2,3,4,5,6\}, \{2,4,5\}, \{3,5,6\}\}$. The bases of the resulting transversal matroid  $M^*$ are all $3$-subsets of $\{1,2,3,4,5,6\}$ except $124$ and $136$.
\end{example}
\medskip

Notice that the cotransversal matroid $M$ of Example \ref{ex1} is dual to the transversal matroid $M^*$ of Example \ref{ex2}. This is a special case of a general phenomenon:

 \begin{theorem}  \label{duals}\cite{ardila, Ingleton, oxley}
Cotransversal matroids are precisely the duals of transversal matroids.
\end{theorem}

Cotransversal matroids were originally called \emph{strict gammoids}. Ingleton and Piff's discovery of Theorem \ref{duals} prompted their newer, widely adopted name.

\section{\textsf{Swapping}}\label{swap}

In this section we introduce the \emph{swap} operation on planted graphs, and show that it preserves the cotransversal matroid. 

In a planted graph, denote the edge from vertex $i$ to vertex $j$ by $e_{ij}$. 

\begin{definition} Let $(G,B)$ be a planted graph, and let $i\notin B,j\in B$ be such that  $e_{ij}\in G$. The \emph{swap operation $\swap(i,j)$} turns $(G,B)$ into the planted graph $(G,B)_{i\rightarrow j}=(G',B')$ by

$\bullet$  
replacing $e_{ij} \in G$ with $e_{ji} \in G'$,

$\bullet$  replacing every other edge of the form $e_{ik}$ in $G$ with $e_{jk} \in G'$, and

$\bullet$ replacing the sink $j \in B$ with the new sink $i \in B'$.

\end{definition}
\begin{figure}[h]
\centering
\includegraphics[scale=.55]{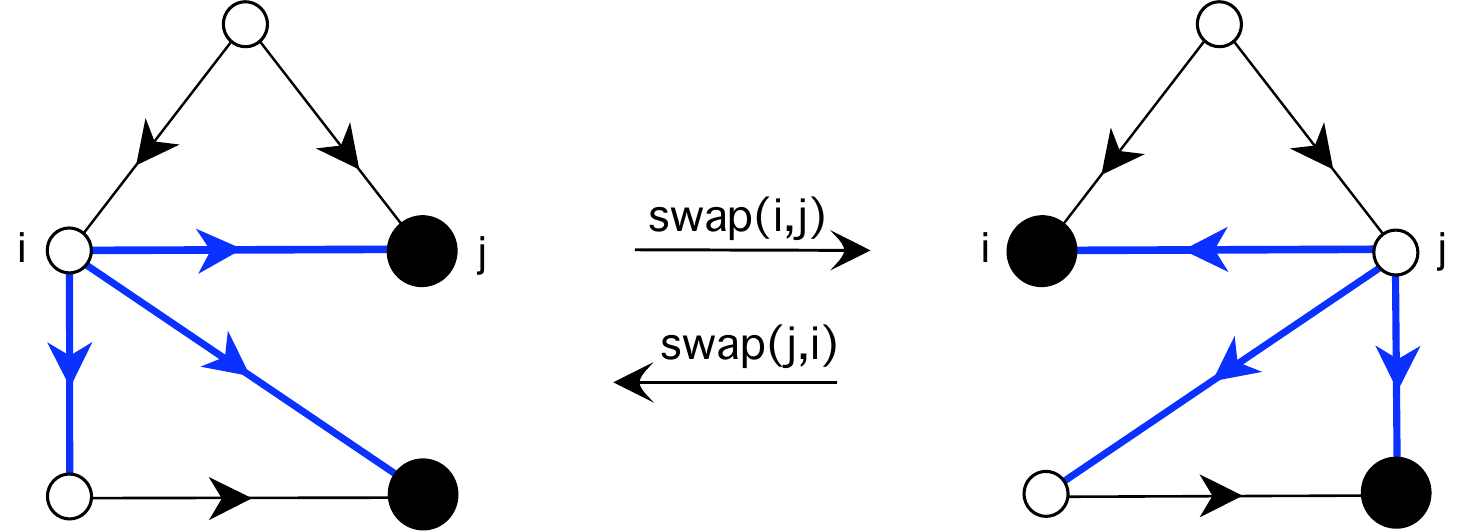} 
\caption{The operation $\swap(i,j)$; sinks are drawn as large black vertices.} \label{swapij}
\end{figure} 

Figure \ref{swapij} illustrates the operation $\swap(i,j)$; the set $B$ is represented by large, black vertices. Notice that $\swap(j,i)$ is a two-sided inverse of $\swap(i,j)$.

\begin{theorem}\label{preserve} 
Swaps preserve the cotransversal matroid: If $(G,B)$ is a planted graph, and $i\notin B,j\in B$ are such that  $e_{ij}\in G$, then $L((G,B)_{i\rightarrow j})=L(G,B)$.
\end{theorem}

\begin{figure}[h]
 \centering
\includegraphics[scale=.65]{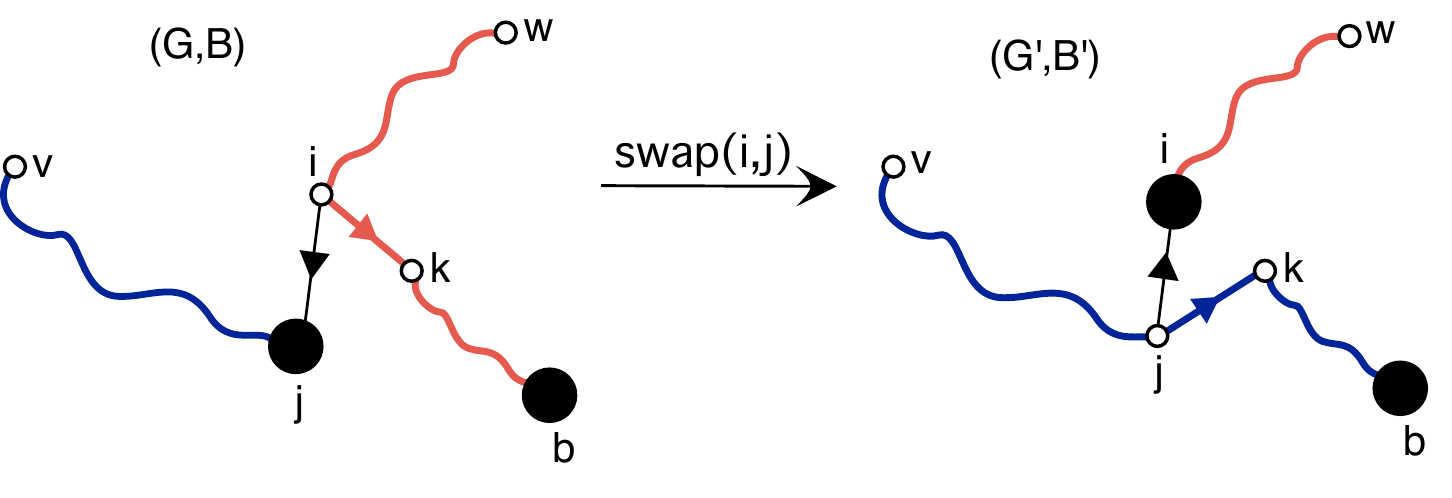}
\caption{Case (iii): Rerouting $v$ and $w$.\label{case3}}
\end{figure}

\begin{proof} 

Since $\swap(i,j)$ is invertible, it suffices to show that any set of vertices which could be routed to $B$ in $(G,B)$ can be routed to $B'$ in $(G,B)_{i \rightarrow j} = (G', B')$.

Let $A$ be a basis of $L(G,B)$, and consider a routing $R$ from $A$ to $B$. Let $p_{ab}$ be the path in $R$ which goes from $a$ to $b$, and let $v$ be the vertex of $A$ which gets routed to $j$. We consider three cases:
(i) $v$ is routed through $i$ to get to $j$,
(ii) $v$ is routed  to $j$ without going through $i$, and $i$ is not in any other route of $R$, and
(iii) $v$ is routed to $j$ without going through $i$, and $i$ is in some other route of $R$.

(i) 
Since $e_{ij}$ is in $G$, we can assume that $R$ uses the path $p_{vj}=(v,\ldots ,i,j)$  from
$v$ to $j$. 
As a result of the operation $\swap(i,j)$ we have $B'=B-j\cup i$. The operation $\swap(i,j)$ does not affect the path from $v$ to $i$, or any other paths in $R$. We can replace the path $p_{vj}$ in $R$ with the path $p'_{vi}=p_{vj}-e_{ij}$ of $G'$, and let the other paths of the routing stay the same.
Therefore $ A$ is a basis of $L(G',B').$

(ii) 
Since $i$ is not on the route from $v$ to $j$, no edges along the
path $p_{vj}$ are affected by the swap, so $v$ still has this path to $j$ in $G'$. Also 
$e_{ji}\in G'$, so the path
 $p'_{vi}=p_{vj} \cup e_{ji}$ in $(G',B')$ routes $v$ to $i$ and doesn't intersect the other paths of the routing. 
 We obtain that $A$ is a basis of $L(G',B').$

(iii) Let $w$ be the vertex of $A$ which is routed through $i$ to some sink $b \in B$, $b \neq j$, as shown in Figure \ref{case3}. 
 As a result of $\swap(i,j)$, the path $p_{wb}$ in $(G,B)$ gets blocked at the edge $e_{ik}$. We can use the truncated path $p'_{wi}=(w,\ldots, i)$ in $(G',B')$ as a route from $w$ to $i \in B'$. 
 To complete a routing we need a path leaving $v \in A$ and arriving at $b \in B'$.
 The path  $p_{vj}$ in $G$
 is unaffected in $G'$, and $e_{jk}\in G'$ since $e_{ik}\in G$. So we can use the old path $p_{vj}$ and the new edge $e_{jk} \in G'$ to pick up the old path from $k$ to $b$; this does not intersect any other path in the routing $R$. 
 It follows that  $A$ is a basis of $L(G',B').$
 \end{proof} 

\section{\textsf{Saturation for cotransversal matroids}}\label{sec:max}

In this section we will see that every presentation $(G,B)$ of a cotransversal matroid $M=L(G,B)$ can be ``saturated" in a unique way into a maximal planted graph $\overline{(G,B)} \supseteq (G,B)$ such that $M=L\overline{(G,B)}$. This is done by adding to $(G,B)$ all missing edges that will not affect the cotransversal matroid. This was essentially proved in \cite{bondy, mason}; to explain it, we need to take a closer look at the duality between cotransversal and transversal matroids.

\subsection{\textsf{Duality between transversal and cotransversal matroids revisited}}

In Theorem \ref{duals} we saw that transversal matroids and cotransversal matroids are dual to each other. We will need a slightly stronger version of this statement:

\begin{theorem}  \label{dualpresentations} \cite{Ingleton}
Let $M$ and $M^*$ be a pair of dual cotransversal and transversal matroids on $V$. Then there is a bijection that maps a planted graph presentation of $M$ to a presentation of $M^*$ together with an SDR. 
\end{theorem}

The previous theorem is implicit in \cite{Ingleton}. For that reason we omit its proof, but we describe the bijection.

Given a planted graph presentation $(G,B)$ of $M$,  let  $A_i:=\{i\}\cup\{u\mid e_{iu}\in G\}$ for each $i \in V-B$. The sets $A_i$ with $i\in V-B$ make up a presentation of $M^*$, and the matching of $i$ with $A_i$ is an SDR for those sets.
 
In the opposite direction, consider a presentation 
$\A=\{A_1,A_2,\ldots , A_k\}$ of $M^*$ and an SDR $a_1, \ldots, a_k$. 
For each $x \in A_j$ with $x \neq a_j$, draw the directed
edge from $a_j$ to $x$ in $G$. Let $B$ be the complement of $\{a_1, \ldots, a_k\}$. This will give a
presentation of $M$.

The reader may find it instructive to check that the planted graph presentation of $M$ in Example \ref{ex1} is dual to the presentation of $M^*$ in Example \ref{ex2} with SDR $(1,2,3)$.

\subsection{\textsf{Saturating a graph}}

As mentioned in Section \ref{sec:prelim}, theorems about a matroid $M$ can often be translated automatically into ``dual" theorems about the dual matroid $M^*$. This is very useful for our purposes. In their foundational work on transversal matroids, Bondy \cite{bondy} and Mason \cite{mason} explained how the different presentations of a transversal matroid are related to each other. Using Theorem \ref{dualpresentations}, we will now ``dualize" their work, to obtain for free several useful results about the presentations of a cotransversal matroid.

The statements in this section are not difficult to show directly. Since they are dual to results in \cite{bondy} and \cite{mason}, we omit their proofs.

\begin{theorem} \label{unique} \cite{bondy, mason} 
For any planted graph $(G,B)$ there exists a unique maximal planted graph $\overline{(G,B)}$ containing $(G,B)$ such that $L\overline{(G,B)} = L(G,B)$. We call $\overline{(G,B)}$ the \emph{saturation} of $(G,B)$.
\end{theorem} 

Theorem \ref{unique} is all that we need to prove our main result, Theorem \ref{iso}. In the rest of this section, which is logically independent from the remainder of the paper, we describe \textbf{how} one constructs the saturation $\overline{(G,B)}$ of $(G,B)$. First we need some definitions.


\begin{definition}
 \label{contract}  
 Let $M=(E, \B)$ be a matroid. Let $K\subseteq E$ and let $B_K$ be a basis of $K$. The \emph{contraction of $M$ by $K$}, denoted $M/K$, is the matroid on $E-K$ whose bases are the sets $B' \subseteq E-K$ such that $B' \cup B_K$ is a basis of $M$.
 \end{definition}
 
It is known \cite[Chapter 5]{white} that any contraction $L(G,B)/K$ of a cotransversal matroid is also cotransversal. To obtain an explicit presentation of it, we first need a presentation $(G',B')$ of $L(G,B)$ with $|K \cap B'| = r(K)$, where $r(K)$ is the maximum number of paths in a routing from $K$ to $B$ in $(G,B)$.
To construct it, start with the planted graph $(G,B)$. If $|K \cap B| < r(K)$, there must be a path from some $k \in K$ to some $b \in B-K$. Performing successive swaps on the edges along this path, one obtains a new presentation $(G_1, B_1)$ where $B_1=B-b \cup k$ satisfies $|K \cap B_1| > |K \cap B|$. By repeating this procedure, we will eventually reach a presentation $(G',B')$ of the matroid with $|K \cap B'| = r(K)$.

%
Finally, delete from $(G',B')$ the vertices in $K$ and all the edges incident to them. It is easy to check that the resulting planted graph is a presentation of the contraction $L(G,B)/K$. 

\begin{definition}Let $v$ be a vertex of a planted graph $(G,B)$. The \emph{claw} of $v$ in $(G,B)$ is $K_v = v\cup\{i\mid e_{vi}\in G \} $. 
\end{definition}

Recall that a \emph{loop} in a matroid is an element that does not occur in any basis of the matroid. In a cotransversal matroid $L(G,B)$, a loop is a vertex of $G$ from which there is no path to $B$. The following proposition tells us which edges we can add to $(G,B)$ without changing the cotransversal matroid.


\begin{proposition} \label{max} \cite{bondy,mason} 
Let $(G,B)$ be a planted graph and let $v$ and $w$ be two vertices of $G$ with $v \notin B$. Then $L(G\cup e_{vw},B) = L(G,B)$ if and only if $w$ is a loop in $L(G,B)/K_v$. 
\end{proposition}

Therefore, to construct the saturation $\overline{(G,B)}$ of a planted graph $(G,B)$, one successively \emph{saturates each vertex} $v \notin B$ as follows: one contracts the matroid by the claw $K_v$, finds the loops in the resulting planted graph, and connects $v$ to those loops. 
In Proposition \ref{max}, the condition for adding the edge $e_{vw}$ depends only on the matroid $L(G,B)$ and the claw $K_v$, neither of which is affected by the saturation of a different vertex $v' \neq v$.  It follows that one can saturate the vertices in any order, and one will always end up with the same graph $\overline{(G,B)}$.

\section{\textsf{An exchange lemma for transversal matroids}}\label{sec:dual}

\begin{theorem} \label{uniquemaxtransv}\cite{bondy, mason}
A transversal matroid has a unique maximal presentation: For every family $\A=\{A_1, \ldots, A_n\}$ of subsets of a set $S$ there is a unique family $\overline{\A} = \{\overline{A_1}, \ldots, \overline{A_n}\}$ of inclusion-maximal subsets of $S$ such that $A_i \subseteq \overline{A_i}$ for $1 \leq i \leq n$, and $\A$ and $\overline{\A}$ give rise to the same transversal matroid.
\end{theorem}

The following lemma on SDRs will be crucial later on.

\begin{lemma}[SDR exchange lemma]\label{match} Suppose that  $\A=\{A_1, \ldots, A_r\}$ satisfies the \emph{dragon marriage condition:}\footnote{This name is due to Postnikov, and originates as follows. Suppose that $S$ is the set of women and $\{1, \ldots, r\}$ is the set of men in a village, and let $A_i$ be the set of women who are willing to marry man $i$. A dragon comes to the village and takes one of the women. When is it the case that all the men can still get married, regardless of which woman the dragon takes away? Postnikov showed that this is the case if and only if $\A$ satisfies the dragon marriage condition.}
 for all nonempty sets $\{i_1, \ldots, i_k\} \subseteq [r]$ we have $|A_{i_1}\cup A_{i_2}\cup\ldots\cup A_{i_k}|\geq k+1$. Then for any two SDRs $M$ and $M'$ of $\A$, there is a sequence $M=M_1, \ldots, M_s=M'$  of SDRs of $\A$ such that $M_i$ and $M_{i+1}$ differ in exactly one position for $1 \leq i \leq s-1$.
\end{lemma}

\begin{proof}  Construct a graph $H$ in which the vertices are the SDRs of $\A$ and two SDRs are connected by an edge if they differ in only one position. We need to prove that $H$ is connected.

Suppose $H$ is not connected.
Consider two SDRs $M_b=(b_1,\ldots,b_r)$ and $M_c=(c_1,\ldots ,c_r)$ in distinct components of $H$. Assume $M_b$ and $M_c$ are chosen so that the Hamming distance $|M_b-M_c|$, \emph{i.e.} the number of positions where $M_b$ and $M_c$ differ, is minimal. We consider the following two cases. 
 
(i) If $\{b_1, \ldots, b_r\} \neq \{c_1, \ldots, c_r\}$, then for some $i$ we have $b_i\notin  \{c_1,\ldots,c_r\}. $ Then $M_c'=(c_1,\ldots,b_i,\ldots,c_r)$ is an SDR in the connected component of $M_c$, and satisfies $|M_b-M'_c| < |M_b-M_c|$.     
 \begin{figure}[h]
 \centering
\includegraphics[scale=.85]{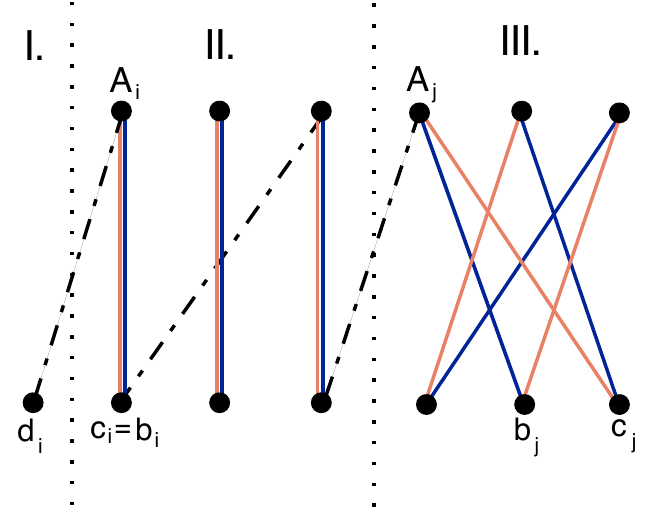}
\caption{Case (ii): $T$ is partitioned into three parts according to the blue and red SDRs $M_b$ and $M_c$.}\label{bip}
\end{figure}

 (ii)  Suppose $\{b_1,\ldots,b_r\}= \{c_1,\ldots,c_r\}$.  
We can partition the vertices of our bipartite graph $T$ into three parts based on the matchings $M_b$ and $M_c$, as shown in Figure \ref{bip}. (The dotted edges will be explained later.)  Part I consists of the vertices of $T$ that are neither in $M_b$ nor in $M_c$. Part II consists of the top vertices $i$ 
 such that $b_i=c_i$, and the bottom vertices matched to them. Part III consists of the remaining vertices.

 The dragon marriage condition gives $|S|\geq r+1$, 
so there is some $d_i\in A_i$ such that $d_i\notin \{b_1,\ldots,b_r\}.$  Therefore $M'_b=(b_1,\ldots,d_i,\ldots,b_r)$ and $M'_c=(c_1,\ldots,d_i,\ldots,c_r)$ are SDRs which are in the connected components of $M_b$ and $M_c$. We must have $b_i=c_i$, or else $|M'_b- M'_c|<|M_b-M_c|$. In Figure \ref{bip}, this means that there are no edges from the top of Part III to Part I.

By the dragon marriage condition, the top of Part III must be connected to the bottom of Part II. Define a \emph{zigzag path} to be a path such that: 

$\bullet$ its starting point is a vertex in the top of Part III, 

$\bullet$ this is the only vertex of Part III it contains, and 

$\bullet$ every second edge is a common edge of the matchings $M_b$ and $M_c$.

We claim that there is at least one zigzag path that ends in Part I.
To verify this, consider the set $U$ of vertices in the top that can be reached by a zigzag path starting from the top of Part III.  Notice that every top vertex in Part III is in $U$.  By the dragon marriage condition, some vertex in $U$ must be connected to a vertex $d$ in the bottom of the graph that is not matched to $U$ in $M_b$ and $M_c$. If $d$ was in Part II,  it would be matched in $M_b$ and $M_c$ to a top vertex $A \notin U$; the edge from $d$ to $A$ would complete a zigzag path that contains $A$, contradicting our definition of the set $U$. Therefore $d$ is in Part I.

Consider a zigzag path to $d$ starting at $A_j$, as shown in Figure \ref{bip}. Now construct new SDRs $M_b'$ and $M_c'$ by unlinking $b_j$ and $c_j$ from $A_j$ in $M_b$ and $M_c$ respectively, as well as all the edges of $M_b$ and $M_c$ along the zigzag path $P$. Instead, in both $M_b$ and $M_c$, rematch the vertices along the edges of path $P$ which were not used by $M_b$ and $M_c$; these are dotted in Figure \ref{bip}. Figure \ref{bip2} shows the resulting new matchings $M_b'$ and $M_c'$ in this example.
\begin{figure}[h]
 \centering
\includegraphics[scale=.85]{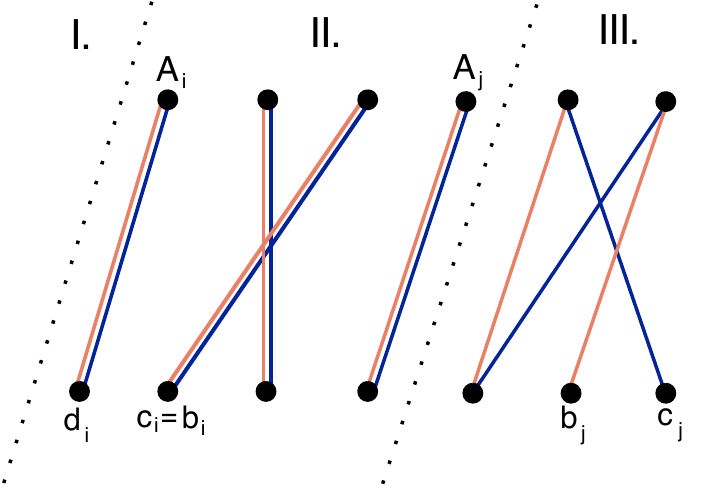}
\caption{The new matchings $M_b'$ and $M_c'$.}\label{bip2}
\end{figure}
Now notice that $|M'_b- M'_c|<|M_b-M_c|$, and $M_b'$ and $M_c'$ are in the same connected components of $H$ as $M_b$ and $M_c$, respectively. This is a contradiction, and we conclude that $H$ is connected. 
\end{proof}


 \section{\textsf{The main result}}\label{proof}

We have now laid all the necessary groundwork to present our main theorem.

\begin{theorem}\label{iso}
Two planted graphs $(G,B)$ and $(H,C)$ have the same cotransversal matroid if and only if  their saturations $\overline{(G,B)}$ and $\overline{(H,C)}$ can be obtained from each other  by a series of swaps. 
\end{theorem}


\begin{proof}
The backward direction follows from Theorems \ref{preserve} and \ref{unique}.
Now suppose $(G,B)$ and $(H,C)$ are presentations of the same cotransversal matroid $M$. When we apply the bijection of Theorem \ref{dualpresentations} to them, both saturations $\overline{(G,B)}$ and  $\overline{(H,C)}$ must give rise to the unique maximal presentation $\A$ of the dual transversal matroid $M^\ast$. They correspond to different matchings $M_1$ and $M_2$ of $\A$. 

Since $\A$ has at least one matching, we have $|A_{i_1} \cup \cdots \cup A_{i_k}| \geq k$ for all $\{i_1, \ldots, i_k\}$ by Hall's theorem. If we have $|A_{i_1} \cup \cdots \cup  A_{i_k}| = k$ for some $\{i_1, \ldots, i_k\}$, then all the elements of $A_{i_1} \cup \cdots \cup A_{i_k}$ are in every basis of $M^*$. Such elements are called \emph{coloops} of $M^*$ and they correspond to loops in $M$. By maximality, the loops of $M$ form a complete subgraph in both $\overline{(G,B)}$ and $\overline{(H,C)}$. This is because loops have no path to the sinks; so they cannot be connected to vertices having paths to the sinks, but they can have any possible connection among themselves. We can then restrict our attention to the non-loops of $M$, where the dragon marriage condition is satisfied.

Applying Lemma \ref{match}, we can get from $M_1$ to $M_2$ by exchanging one element of the matching at a time. One easily checks that these matching exchanges in the bipartite graph correspond exactly to swaps in the corresponding planted graphs under the bijection of Theorem \ref{dualpresentations}. It follows that one can get from $\overline{(G,B)}$ to $\overline{(H,C)}$ by a series of swaps, as desired. 
\end{proof}

We end by illustrating Theorem \ref{iso} with two examples.

\begin{figure}[h]
\centering
\includegraphics[scale=.5]{mygammoid_123max}
\hspace{1cm}\includegraphics[scale=.5]{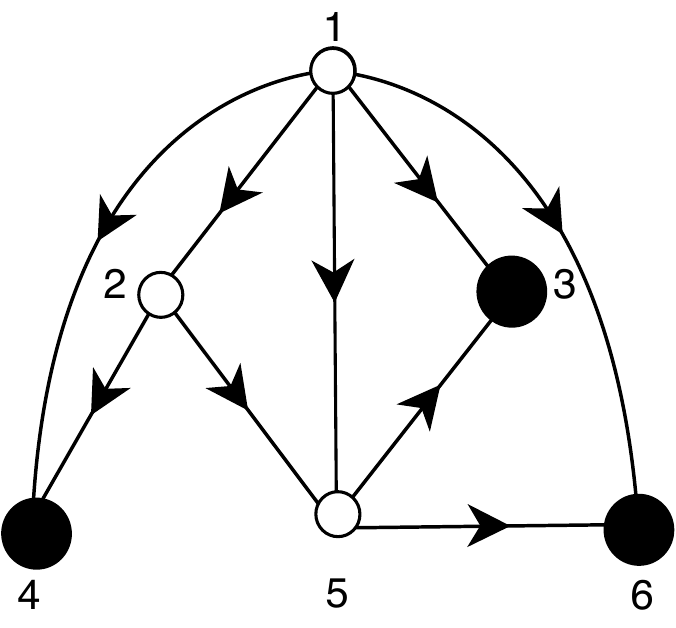}
\hspace{1cm}\includegraphics[scale=.5]{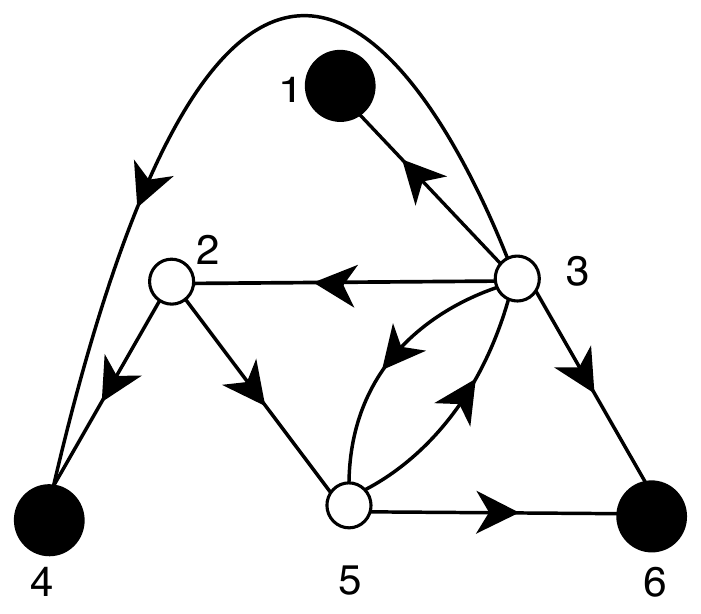}
\caption{
The planted graphs given by 
$\A = \{\{1, 2, 3, 4, 5, 6\}, \{2, 4, 5\}, \{3, 5, 6\}\}$ with SDRs $(1, 2, 3), (1,2,5)$, and $(3, 2, 5)$, respectively.
\label{moves}}
\end{figure}

\begin{example}
Figure \ref{moves} shows three saturated planted graph presentations of the cotransversal matroid of Example \ref{ex1}. They correspond to the dual maximal presentation $\A = \{\{1, 2, 3, 4, 5, 6\}, \{2, 4, 5\}, \{3, 5, 6\}\}$  of the transversal matroid of Example \ref{ex2}, with SDRs $(1, 2, 3), (1,2,5)$, and $(3, 2, 5)$, respectively. Notice how one-position exchanges in the SDRs correspond to swaps in the planted graphs.
\end{example}

\begin{example}
Let $M$ be the cotransversal matroid on $\{1,2,3,4,5\}$ with bases $\{14, 15, 24, 25, 34, 35, 45\}$. Figure \ref{fig:graph} shows the graph of saturated planted graph presentations of $M$, where two planted graphs are joined by an edge labelled $ij$ if they can be obtained from one another by $\swap(i,j)$. There are nine saturated presentations in two isomorphism classes. We have drawn one representative from each isomorphism class; every other saturated presentation is obtained from one of these two planted graphs by relabelling the vertices. 

\begin{figure}[h]
\centering
\includegraphics[scale=.55]{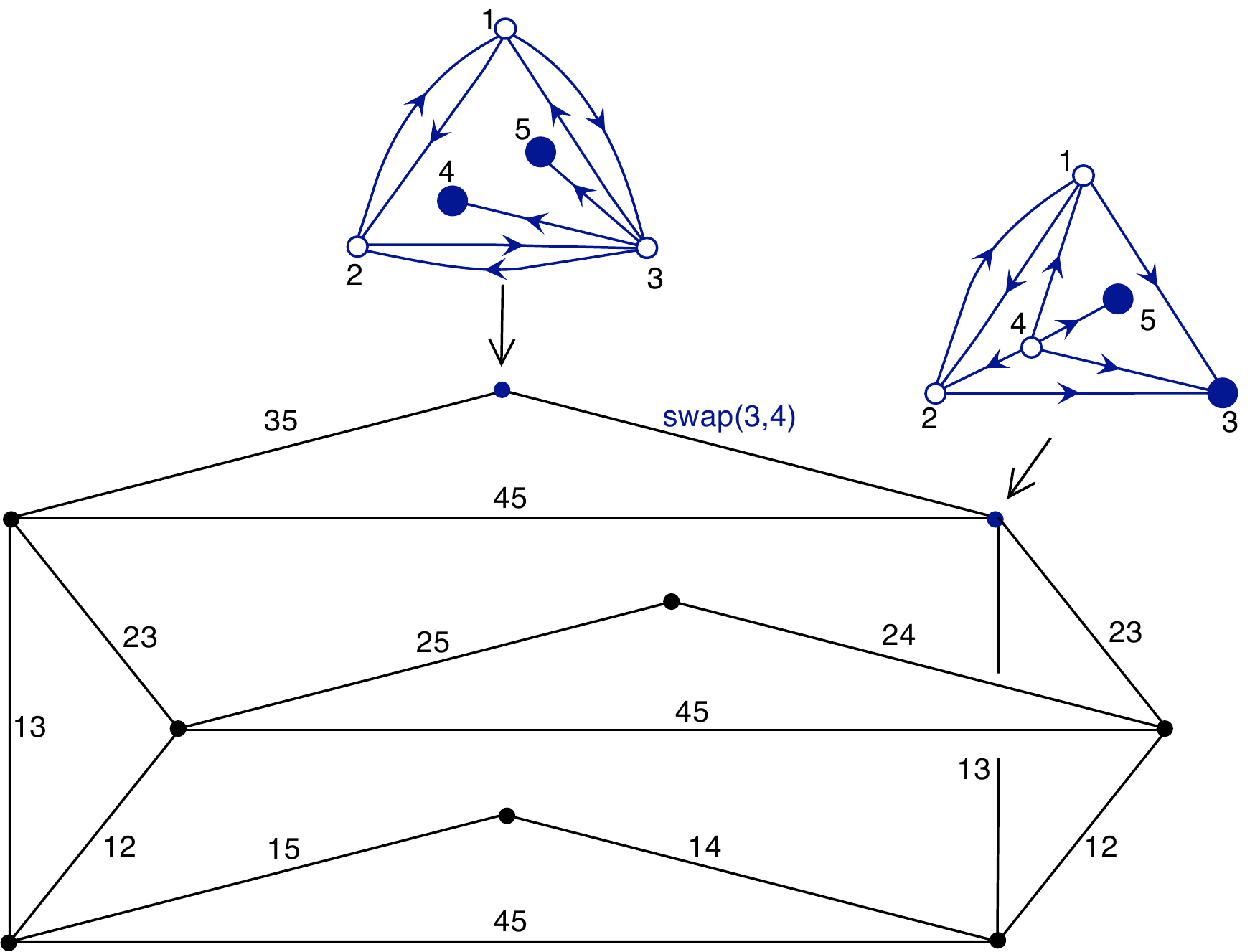}
\caption{The graph of saturated presentations of a cotransversal matroid.}\label{fig:graph}
\end{figure}
\end{example}

\section{\textsf{Acknowledgments}}
We would like to thank the referee for a thorough report and helpful suggestions to improve the presentation.

\footnotesize{

}

\nocite{*}
\begin{thebibliography}{1}
\providecommand{\natexlab}[1]{#1}
\providecommand{\url}[1]{\texttt{#1}}
\expandafter\ifx\csname urlstyle\endcsname\relax
  \providecommand{\doi}[1]{doi: #1}\else
  \providecommand{\doi}{doi: \begingroup \urlstyle{rm}\Url}\fi

  
\bibitem{ardila} F.~Ardila. Transversal and cotransversal matroids via their representations. \emph{Electronic Journal of Combinatorics} {\bf14} (2007), \#N6.


\bibitem{ardilavideos} F.~Ardila. SFSU-Los Andes lecture notes and videos on Matroid Theory, 2007. Available at \textsf{http://math.sfsu.edu/federico/matroids.html}.

\bibitem{bondy}J.A.~Bondy. {Presentations of transversal matroids}. \emph{Journal of the London Mathematical Society} {\bf 5} (1972), 289-292.



\bibitem{Ingleton}
A.~Ingleton and M.~Piff. Gammoids and transversal matroids.
\emph{J. Combinatorial Theory Ser. B} {\bf 15} (1973) 51-68.

\bibitem{mason}J. H. Mason. {Representations of independence spaces}. \emph{University of Wisconsin PhD thesis.} (1970).

\bibitem{mason2} J. H. Mason. On a class of matroids arising from paths in graphs. \emph{Proc. London Math. Soc.} (3) {\bf 25}
(1972) 55-74.

\bibitem{oxley} J.~Oxley. \emph{Matroid Theory}. Oxford University
Press. New York, 1992.


\bibitem{Postnikov}
A. Postnikov, Permutohedra, associahedra, and beyond. \texttt{arXiv:math.CO/0507163}. 
To appear in \emph{International Mathematics Research Notices.} 


\bibitem{white}
N. White (ed.). \emph{Combinatorial geometries}. Cambridge University Press. Cambridge, 1987.

\bibitem{whitney}
H. Whitney. 2-isomorphic graphs. \emph{Amer. J. Math.} {\bf 55} (1933), 245-254.
\end{thebibliography}
\end{document}